\documentclass[review]{elsart}

\usepackage{natbib, amssymb, latexsym, amsmath}
\begin{document}
\begin{frontmatter}
\title{Weight Ideals Associated to Regular and Log-Linear Arrays}
\author{Jeremiah W. Johnson}
\address{Dept. of Mathematics and Statistics, Penn State Harrisburg \\
Middletown, PA 17057}
\ead{jwj10@psu.edu}

\begin{abstract}
Certain weight-based orders on the free associative algebra $R = k\langle x_1, \dots, x_t \rangle$ can be specified by $t \times \infty$ arrays whose entries come from the subring of nonnegative elements in a totally ordered field. Such an array $A$ satisfying certain additional conditions produces a partial order on $R$ which is an admissible order on the quotient $R/I_A$, where $I_A$ is a homogeneous binomial ideal called the {\em weight ideal\/} associated to the array and whose structure is determined entirely by $A$. This article discusses the structure of the weight ideals associated to two distinct sets of arrays whose elements define admissible orders on the associated quotient algebra.
\end{abstract}

\begin{keyword}Noncommutative Gr\"obner Bases, Gr\"obner Bases, Admissible Orders
\end{keyword}
\end{frontmatter}

\section{Introduction} Work over the past two decades has extended the theory of Gr\"obner bases to various noncommutative algebras \citep{EG00, MR97, PN01, TM94}. Before a Gr\"obner basis for an ideal of a $k$-algebra $\mathcal{A}$ can be constructed, where $k$ is a field, an admissible order on a multiplicative basis of $\mathcal{A}$ is required. Following \citet{EG96}, we say that $\mathcal{A}$ has a {\em Gr\"obner basis theory\/} when an admissible order exists on a multiplicative basis of $\mathcal{A}$. In \citet{EH10}, E. Hinson adapted the theory of position-dependent weighted orders to define a length-dominant partial order on the set of words in the free associative algebra $R = k\langle x_1, \dots, x_t \rangle$, including the trivial word, which produces an admissible order on a quotient of $R$. In this construction, the partial order on $R$ is specified by a $t\times\infty$ array $A$ whose entries come from the subring consisting of the positive elements of a totally ordered field, and the quotient is by a homogenous binomial ideal $I_A$ whose elements are determined by the partial order given by $A$. This gives rise to two immediate questions. First, given an array $A$ that defines an admissible order on a quotient $R/I_A$, what is the algebra that is determined, or more specifically, what is the structure of the ideal $I_A$? Second, given two arrays $A$ and $B$ which define orders $\succ_A$ and $\succ_B$ on $R/I_A$ and $R/I_B$ respectively, even when $R/I_A = R/I_B$ it is not necessarily the case that $\succ_A = \succ_B$. Under what circumstances does $\succ_A = \succ_B$? This paper describes results concerning the first of these two questions for two distinct families of admissible arrays. In this introductory section, we review the relevant definitions and results from \citet{EH10} and we make the preceding general statements precise. Our primary objects of interest are defined in Definitions \ref{DEF:WeightIdeal} and \ref{DEF:AdmisArray}. The results on which the remainder of the paper relies are given in Theorems \ref{THM:AdmisArray} and \ref{THM:AdmisOrder}. In what follows, let $R = k\langle x_1, \dots, x_t \rangle$ denote the free associative algebra, let $S_{>0}$ denote the positive elements of a totally ordered field, and let $\mathcal{M}_{t\times \infty}\left(S_{>0}\right)$ denote the set of $t\times\infty$ arrays with entries in $S_{>0}$. The following two definitions are adopted from \citet{EG96}.

\begin{defn}\label{DEF:MultBasis} Let $\mathcal{B}$ be a $k$--basis of an algebra $\mathcal{A}$. $\mathcal{B}$ is a {\em multiplicative basis\/} for $\mathcal{A}$ if
\[
b,b' \in \mathcal{B} \Rightarrow b\cdot b' \in \mathcal{B}\text{ or }b\cdot b'=0.
\] 
\end{defn}

We will have occasion to refer to the nontrivial elements of $\mathcal{B}$, which we denote by $\mathcal{B}^{\times}$.

\begin{defn}\label{DEF:AdmissOrder}
A total order $\succ$ on a multiplicative basis $\mathcal{B}$ of $\mathcal{A}$ is an {\em admissible order on\/} $\mathcal{B}$ if  
\begin{itemize}
\item $\succ$ is a well-order on $\mathcal{B}$,
\item for all $b_1$, $b_2$, $b_3 \in \mathcal{B}$ such that $b_1b_3 \neq 0$ and $b_2b_3 \neq 0$, if $b_1 \succ b_2$, then $b_1b_3 \succ b_2b_3$,
\item for all $b_1$, $b_2$, $b_3 \in \mathcal{B}$ such that $b_3b_1 \neq 0$ and $b_3b_2 \neq 0$, if $b_1 \succ b_2$, then $b_3b_1 \succ b_3b_2$, and
\item for all $b_1$, $b_2$, $b_3$, $b_4 \in \mathcal{B}$, if $b_1 = b_2b_3b_4$, then $b_1 \succeq b_3$.
\end{itemize}
\end{defn}

Commonly used admissible orders for Gr\"obner basis calculations on noncommutative algebras are the left length-lexicographic order or the right length-lexicographic order \citep{EG96}. We specify a position-dependent weighted order on words in the free algebra using a $t\times\infty$ array to define a weight function as described in the following definition.

\begin{defn} Let $A=(a_{i,j}) \in\mathcal{M}_{t\times\infty}(S_{>0})$. $A$ gives a monomial weighting $\sigma_A: \mathcal{B}^{\times} \to S_{>0}$ by
\[
\sigma_A(x_{u_0}x_{u_1}\cdots x_{u_{l-1}}) = \prod_{j = 0}^{l-1}a_{u_{j},j}
\]
    
for a given monomial $x_{u_0}x_{u_1}\cdots x_{u_{l-1}} \in R$. The function $\sigma_A$ is the {\em weight function associated to\/} $A$.
\end{defn}

Note that for computational convenience we index the columns of an array starting with 0 rather than 1. When the array $A$ is clear, we will suppress it from the notation and write the associated weight function $\sigma_A$ simply as $\sigma$. For the remainder of this section, fix an array  $A\in\mathcal{M}_{t\times\infty}(S_{>0})$ and associated weight function $\sigma$. In order to discuss the weight of the product of two words, we identify a translated version of the weight function associated to $A$ by
\[
\sigma_{k}(x_{u_0}x_{u_1}\cdots x_{u_{l-1}}) = \prod_{j = 0}^{l-1}a_{u_{j},j+k},
\]

where $k \in \Nset$. We consider $\sigma(\omega) = \sigma_{A}(\omega) = \sigma_{A,0}(\omega)$. Let $|\omega|$ denote the length of $\omega$. Given $\omega $ and $\lambda$ such that $| \omega |=k$, 
\[
\sigma(\omega\lambda) = \sigma(\omega)\cdot \sigma_k(\lambda).
\]

This gives rise to the following equivalence relation.

\begin{defn}\label{DEF:EquivRel} Define the relation $\succ_{\sigma}$ on $\mathcal{B}^{\times}$ by
\[
\omega_1 \succ_{\sigma} \omega_2 \iff |\omega_1| > |\omega_2| \text{, or } |\omega_1| = |\omega_2| \text{ and } \sigma(\omega_1) > \sigma(\omega_2).
\]
\end{defn}
 
Let $\Gamma$ denote the set of pure homogeneous binomial differences $\omega_1-\omega_2$, where $\omega_1$, $\omega_2\in\mathcal{B}^{\times}$, $|\omega_1| = |\omega_2|$, and $\sigma(\omega_1) = \sigma(\omega_2)$. 

\begin{defn}\label{DEF:WeightIdeal}
The ideal $I_A = \langle \Gamma \rangle$ is the {\em weight ideal associated to\/} $A$.
\end{defn}
  
\begin{defn}\label{DEF:AdmisArray} $A$ is an {\em admissible\/} array if for every pair $\omega_1$, $\omega_2 \in \mathcal{B}^{\times}$ with $|\omega_1| = |\omega_2|$,
\begin{itemize}
\item[(1)] for all $k \geq 0$, if $\sigma_k(\omega_1) > \sigma_k(\omega_2)$, then $\sigma_{k+1}(\omega_1) > \sigma_{k+1}(\omega_2)$, and
\item[(2)] for all $k \geq 0$, if $\sigma_k(\omega_1) = \sigma_k(\omega_2)$, then $\sigma_{k+1}(\omega_1) = \sigma_{k+1}(\omega_2)$.
\end{itemize} 
\end{defn}

The following theorem illustrates that the second part of Definition \ref{DEF:AdmisArray} is in fact unnecessary.

\begin{thm}\label{THM:AdmisArray} Let $A \in \mathcal{M}_{t \times \infty}\left(S_{>0}\right)$ be an array with associated weight function $\sigma$. The following are equivalent:
\begin{itemize}
\item[(1)] $A$ is an admissible array;
\item[(2)] for all $k \geq 0$ and for all $\omega_1$, $\omega_2 \in \mathcal{B}^{\times}$ such that $|\omega_1| = |\omega_2|$, $\sigma_k(\omega_1) > \sigma_k(\omega_2)$ if and only if $\sigma_{k+1}(\omega_1) > \sigma_{k+1}(\omega_2)$.
\end{itemize}
\end{thm}

Admissible arrays define an admissible order on the quotient $R/I_A$.

\begin{thm}\label{THM:AdmisOrder} An array $A\in \mathcal{M}_{t \times \infty}\left(S_{>0}\right)$ with associated weight function $\sigma$ is an admissible array if and only if $\succ_{\sigma}$ is an admissible order on $\mathcal{B}_{\sigma} \subseteq R/I_A$, where $\mathcal{B}_{\sigma}$ is the image of $\mathcal{B}$ in $R/I_A$ under the projection $R \to R/I_A$.
 \end{thm}  
 
  \begin{defn}\label{DEF:Degenerate} $A$ is said to be {\em degenerate\/} if there exists $i$, $j$, $1 \leq i \neq j \leq t$, such that $\sigma(x_i) = \sigma(x_j)$.
 \end{defn}
 
 We will assume in what follows that all arrays considered are nondegenerate, for if $\sigma(x_i) = \sigma(x_j)$ for some $i,j\in\{1,\dots,t\}$ where $i \neq j$, then $x_i -x_j \in I_A$ and $k\langle x_1,\dots,x_t\rangle/I_A \simeq k\langle x_1,\dots,x_{i-1},x_{i+1},\dots,x_t\rangle/\langle I_A' \rangle$ where $A'$ is the array obtained from $A$ by deleting the $i^{th}$ row.

\section{Weight Ideals Associated to Regular Arrays}

In \citet{EH10}, E. Hinson described two sets of admissible arrays. We begin by studying the first of these, the set of regular arrays. 

\begin{defn} An array $A$ is {\em regular\/} if $A$ has rank 1.
\end{defn}

The set of linear arrays is a subset of the set of regular arrays which will be used later on to construct the set of log-linear arrays.  

\begin{defn} An array $A$ is {\em linear\/} if for all $i \geq 1$, $A_{(i)} = d\cdot A_{(i-1)}$ for some fixed $d \in S_{>0}$. The fixed scalar $d$ is referred to as the {\em slope\/} of the array. 
\end{defn}

\begin{exmp} The array 
\[
A=\left(\begin{array}{cccc}
2 & 6 & 18 &\cdots  \\
3 & 9 & 27 & \cdots \\
4 & 12 & 36 & \cdots 
\end{array}\right)
\]
is a linear array with slope $d=3$.
\end{exmp}

The weight ideal associated to a regular array contains the commutator ideal $\mathcal{C} = \langle x_ix_j -x_jx_i | 1 \leq i \neq j \leq t \rangle$, and thus is never trivial \citep{EH10}. 
 
\begin{defn}\label{DEF:Support} The {\em support\/} of a word $\omega$ is the set
\[
	\text{supp}(\omega) = \{x_i| i\in \{1,\dots, t\} \text{ and } x_i \text{ occurs in }\omega\}.
\]
\end{defn}  

\begin{defn}\label{DEF:Frequency} The {\em frequency of $x$ in $\omega$\/} is the number of times that $x$ occurs in $\omega$ and is written $\#(x,\omega)$.
\end{defn}

\begin{defn}\label{DEF:AlgConsequence} Let $f \in R$ and $G = \{g_1,g_2,\dots\} \subset R$. We say that $f$ is an {\em algebraic consequence of\/} $G$ if $f = \sum_{g \in G}c_iu_ig_iv_i$, where $c_i \in k$, $u_i$, $v_i \in R$, and only finitely many $c_i \neq 0$. 
\end{defn}

Suppose $A$ is a regular array with first column $[a_{1,0},\dots, a_{t,0}]^T$, where $a_{i,0} \in \Nset$ and at least one of $(a_{i,0}, a_{j,0}) \neq 1$, where $(a_{i,0}, a_{j,0})$ denotes the greatest common divisor of $a_{i,0}$ and $a_{j,0}$ and $1 \leq i \neq j \leq t$. Let $\omega_1 = x_{u_0}\cdots x_{u_{l-1}}$ and $\omega_2 = x_{v_0}\dots x_{v_{l-1}} \in \mathcal{B}$ such that $\omega_1 -\omega_2 \in I_A$. Then we have
\[
\prod_{i = 0}^{l-1} a_{{u_i}i} = \prod_{i=0}^{l-1}a_{{v_i}i},
\]

and each $a_{{u_i}i}$ and $a_{{v_i}i}$ can be written as scalar multiples of $a_{{u_i}0}$ and $a_{{v_i}0}$ respectively: 
\[
\prod_{i = 0}^{l-1} d_ia_{{u_i}0} = \prod_{i=0}^{l-1}d_ia_{{v_i}0}.
\]

Factoring out and canceling the common $d_i$'s reduces the equation to 
\begin{equation}\label{EQ:LinearEq}
\prod_{i = 0}^{l-1} a_{{u_i}0} = \prod_{i=0}^{l-1}a_{{v_i}0}.
\end{equation} 

Equation (\ref{EQ:LinearEq}) does not depend on the how the variables were ordered in $\omega_1$ and $\omega_2$; in particular, by factoring out and canceling any terms $a_{{u_i}0} = a_{{v_j}0}$ common to both sides of the equation, one obtains the reduced expression
\begin{equation}
\prod_{k = 0}^{n}a_{{u_k}0} = \prod_{k=0}^{n}a_{{v_k}0}.
\end{equation}

In this expression, $a_{{u_m}0} \neq a_{{v_{m'}}0}$ for all $u_m$ and $v_{m'}$. Note that we have not cancelled any common divisors of the $a_{u_i,0}$, we have only cancelled those $a_{u_i,0}$'s and $a_{v_j,0}$'s for which $a_{u_i,0} = a_{v_j,0}$. Since each $a_{{u_m}0}$ and $a_{{v_{m'}}0}$ corresponds to the weight assigned to an individual letter in $\{x_1,\dots,x_t\}$, this equation describes a homogeneous binomial difference $\omega_1' - \omega_2' \in I_A$ in which no letter that occurs in $\omega_1'$ will occur in $\omega_2'$. 

\begin{defn}\label{DEF:DisjointSupp} A homogeneous binomial difference $\omega_1 - \omega_2 \in I_A$ for which
\[
\text{supp}(\omega_1) \bigcap \text{supp}(\omega_2) = \emptyset
\]
will be referred to as a {\em homogeneous binomial difference of disjoint support\/}.
\end{defn}

Homogeneous binomial differences of disjoint support may arise as algebraic consequences of other homogeneous binomial differences of disjoint support. For example, suppose 
 \[
 A=\left(\begin{array}{ccccc}
 2 & 4 & 8 & \cdots \\
 3 & 6 & 12 & \cdots \\
 4 & 8 & 16 & \cdots \\
 6 & 12 & 24 & \cdots 
 \end{array}\right).
  \]
 This array $A$ is linear with slope 2. Consider the homogeneous binomial difference $x_3x_2x_3x_2 - x_4x_1x_4x_1$. Since 
 \[\sigma(x_3x_2x_3x_2) = \sigma(x_4x_1x_4x_1) = 9216,
 \]
 we must have $x_3x_2x_3x_2 - x_4x_1x_4x_1 \in I_A$. Neither word in this homogeneous difference shares a letter with the other, so $x_3x_2x_3x_2 - x_4x_1x_4x_1$ is a homogeneous binomial difference of disjoint support. Furthermore,
 \[
 x_3x_2x_3x_2 - x_4x_1x_4x_1 = (x_3x_2-x_4x_1)x_3x_2 +x_4x_1(x_3x_2-x_4x_1),
 \]
  
so $x_3x_2x_3x_2 - x_4x_1x_4x_1$ is a homogeneous binomial difference of disjoint support which arises as an algebraic consequence of a homogeneous binomial difference of disjoint support consisting of words of lesser length. 

 \begin{defn}\label{DEF:MinDisjoint} Let $\omega_1-\omega_2$ be a homogeneous binomial difference of disjoint support. $\omega_1-\omega_2$ is {\em minimal\/} if any expression 
 \[
  \omega_1-\omega_2 = \sum_{i = 1}^{n}\alpha_i(u_i-v_i)\beta_i
  \]
for $\omega_1-\omega_2$ as a sum of homogeneous binomial differences has at least one difference $u_i-v_i$ such that $|u_i| = |v_i| =|w_1|$. $\mathcal{M}_A$ will be used to denote the set of minimal length homogeneous binomial differences of disjoint support associated to $A$.
 \end{defn}
 
 In other words, a minimal homogeneous binomial difference of disjoint support is one which cannot be realized as an algebraic consequence of homogeneous binomial differences of disjoint support consisting of words of lesser length. 

Any element of $I_A$ may be decomposed over the set of commutators $\{x_ix_j - x_jx_i | 1 \leq i \neq j \leq t\}$ and the set of homogeneous binomial differences of disjoint support.  

\begin{lem}\label{LEM:Algorithm} Let $\omega_1 - \omega_2$ be a homogeneous binomial difference in $I_A$. Then $\omega_1-\omega_2 = \sum_{i = 1}^{n}\alpha_i(u_i-v_i)\beta_i$, where each homogeneous binomial difference $u_i-v_i$, $1\leq i < n$ is a commutator and $u_n-v_n$ is a homogeneous binomial difference of disjoint support.
\end{lem}

\begin{pf} Suppose $\omega_1- \omega_2 \in I_A$. We proceed by induction. The base case when $l = 2$ is established trivially. Assume now that the induction hypothesis holds for homogeneous binomial differences consisting of words of length $l-1$ and suppose $|\omega_1|=|\omega_2| = l$. Write $\omega_1=x_{u_0}\dots x_{u_{l-1}}$ and $\omega_2=x_{v_0}\dots x_{v_{l-1}}$. Let $i \in \{0,\dots, l-1\}$ be the least value for which $x_{u_0} = x_{v_i}$ (if no such value exists, we are done). By inserting the expression
\[
-x_{v_0}\dots x_{v_{i-2}}x_{v_i}x_{v_{i-1}}x_{v_{i+1}}\dots x_{v_{l-1}}+x_{v_0}\dots x_{v_{i-2}}x_{v_i}x_{v_{i-1}}x_{v_{i+1}}\dots x_{v_{l-1}},
\]

we obtain
\[ 
\omega_1-x_{v_0}\dots x_{v_{i-2}}x_{v_i}x_{v_{i-1}}x_{v_{i+1}}\dots x_{v_{l-1}}+x_{v_0}\dots x_{v_{i-2}}x_{v_i}x_{v_{i-1}}x_{v_{i+1}}\dots x_{v_{l-1}}-\omega_2,
\]

which is equal to
\begin{multline}\label{EQ:FirstFact}
\omega_1-x_{v_0}\dots x_{v_{i-2}}x_{v_i}x_{v_{i-1}}x_{v_{i+1}}\dots x_{v_{l-1}}+ \\
	x_{v_0}\dots x_{v_{i-2}}\left(x_{v_i}x_{v_{i-1}}-x_{v_{i-1}}x_{v_i}\right)x_{v_{i+1}}\dots x_{v_{l-1}}.
\end{multline}

In the second term in expression (\ref{EQ:FirstFact}), $x_{v_i}$ occurs in the $i-1^{st}$ position. The third and fourth terms in Equation \ref{EQ:FirstFact} have been expressed as (left and right) multiples of the commutator $x_{v_i}x_{v_{i-1}}-x_{v_{i-1}}x_{v_i}$. Iterating this process $i$ times results in the expression
\begin{equation}\label{EQ:SecondFact}
\omega_1-x_{v_i}x_{v_0}\dots x_{v_{i-1}}x_{v_{i+1}}\dots x_{v_{l-1}}+\sum_{k=1}^{i-1}\alpha_k(x_{v_i}x_{v_{i-k}}-x_{v_{i-k}}x_{v_i})\beta_k,
\end{equation}

where $\alpha_k = x_{v_0} \dots x_{v_{i-k-1}}$ and $\beta_k = x_{v_{i-k+1}} \dots x_{v_{l-1}}$.

Since $x_{u_0} = x_{v_i}$, the difference of the first two terms in \ref{EQ:SecondFact} can be rewritten as 
\[
x_{u_0}\left(x_{u_1}\dots x_{u_{l-1}}-x_{v_1}\dots x_{v_{l-1}}\right).
\]

The expression in parentheses consists of monomials of length $l-1$ which is an algebraic consequence of the commutators and a homogeneous binomial difference of disjoint support. Rearranging and renaming terms as needed gives the desired result. \qed 
\end{pf}
 
 \begin{thm}\label{THM:RegGen} Let $A$ be a regular array. The weight ideal $I_A$ associated to a regular array $A$ is generated by the union of the set of commutators $\{x_ix_j-x_jx_i| 1 \leq i \neq j \leq t\}$ and $\mathcal{M}_A$.
 \end{thm} 
 
 \begin{pf} Fix a homogeneous binomial difference $\omega_1-\omega_2 \in I_A$. By iterating the algorithm described in the proof of Lemma \ref{LEM:Algorithm}, $\omega_1-\omega_2 \in I_A$ can be reduced to an algebraic consequence of the commutators plus a single, perhaps trivial, homogeneous binomial difference of disjoint support $\omega_1' -\omega_2'$. To see this, note that each iteration of the algorithm produces in the sum a difference of commutators and a homogeneous binomial difference of shorter length than in the previous iteration in which a letter common to each word has been extracted. We may continue the algorithm until either the next iteration is over a commutator or there are no common letters to extract. In the first case, we are done, and in the second case, if $\omega_1' -\omega_2'$ is minimal, we are also done. If $\omega_1' -\omega_2'$ is not a minimal homogeneous binomial difference, then by definition it is an algebraic consequence of minimal homogeneous binomial differences of disjoint support. \qed
 \end{pf}

Having obtained a description of the generators of $I_A$, we will next show that when $A$ is regular, $I_A$ is finitely generated. We include the following lemma to describe the means by which a disjoint homogeneous binomial difference which contains another difference as scattered subwords can be decomposed over that subdifference.

\begin{lem}\label{LEM:Subwords} Let $\omega_1 - \omega_2 \in \mathcal{M}_A$ and suppose $\lambda_1 -\lambda_2$ is a homogeneous binomial difference of disjoint support such that $\omega_1$ occurs as a scattered subword in $\lambda_1$ and $\omega_2$ occurs as a scattered subword in $\lambda_2$. Then 
\[	
	\lambda_1 - \lambda_2 =  (\omega_1 - \omega_2)\alpha+\omega_2(\alpha - \beta) + \sum_{i = 1}^n \alpha_i(\gamma_i - \zeta_i)\beta_i,
\]
where $\alpha - \beta$ is a homogeneous binomial difference of disjoint support and $\gamma_i - \zeta_i$ is a commutator for each $i$, $1\leq i \leq n$.
\end{lem}

\begin{pf} The algorithm of Lemma \ref{LEM:Algorithm} may be modified to move any letter that occurs in a word in a homogeneous binomial difference in $I_A$ either forward or backwards to the desired position, resulting in a decomposition
\[
	\lambda_1 -\lambda_2 = \omega_1\alpha - \omega_2\beta + \sum_{i = 1}^n \alpha_i(\gamma_i - \zeta_i)\beta_i,
\]
where $\gamma_i - \zeta_i$ is a commutator, $1\leq i \leq n$. The result then follows.	
 \end{pf}
 
\begin{thm}\label{THM:FinGen} Let $A$ be a regular array. The associated weight ideal $I_A$ is finitely generated.
\end{thm}
  
\begin{pf} By Theorem \ref{THM:RegGen}, $I_A$ is generated by the union of the set of commutators and the set $\mathcal{M}_A$ of minimal homogeneous binomial differences of disjoint support. The set of commutators is clearly finite. It remains to demonstrate that $\mathcal{M}_A$ is also finite. Assume the contrary. Then there exists some partition of $X = \{x_1,\dots,x_t\}$ into two sets $X_1$, $X_2$ such that there are infinitely many minimal disjoint homogeneous binomial differences $\omega_1 - \omega_2$ in which $\text{supp}(\omega_1) \subseteq X_1$ and $\text{supp}(\omega_2) \subseteq X_2$. Let $\mathcal{D} = \{\omega_1-\omega_2\in \mathcal{M}_A| \text{supp}(\omega_1) \in X_1\text{, }\text{supp}(\omega_2) \in X_2\}$ and let $\omega_1 - \omega_2 \in \mathcal{D}$ such that $|\omega_1| \leq |\lambda_1|$ for any $\lambda_1$ that occurs in a homogeneous binomial difference $\lambda_1 - \lambda_2 \in \mathcal{D}$. Consider the following three sets: $\mathcal{D}(\omega_1) = \{\lambda_1 - \lambda_2 \in \mathcal{D} : \omega_1\text{ occurs as a scattered subword in }\lambda_1\}$, $\mathcal{D}(\omega_2) = \{\lambda_1 - \lambda_2 \in\mathcal{D} : \omega_2\text{ occurs as a scattered subword in } \lambda_2\}$, and $\mathcal{D}(0) = \{\lambda_1-\lambda_2\in\mathcal{D} : \text{neither }\omega_1 \text{ nor } \omega_2 \text{ occur as scattered subwords in } \lambda_1\text{ and }\lambda_2\}$. Note that $\mathcal{D} = \{\omega_1-\omega_2\}\cup\mathcal{D}(\omega_1)\cup\mathcal{D}(\omega_2)\cup\mathcal{D}(0)$. Furthermore, these sets are disjoint. If $\omega_1$ were to occur as a scattered subword in $\lambda_1$ and $\omega_2$ occurs as a scattered subword of $\lambda_2$, then Lemma \ref{LEM:Subwords} shows that $\lambda_1 - \lambda_2$ is an algebraic consequence of commutators, $\omega_1-\omega_2$, and perhaps some other homogeneous binomial difference in $\mathcal{D}$ consisting of words of length less than $|\lambda_1|$; that is, $\lambda_1 - \lambda_2$ is not minimal. Thus, these sets form a partition of $\mathcal{D}$ and so at least one of $\mathcal{D}(\omega_1)$, $\mathcal{D}(\omega_2)$, and $\mathcal{D}(0)$ must be infinite.

Now let $\lambda_1 -\lambda_2 \in \mathcal{D}(\omega_2)$ and suppose $|\lambda_1| > |\omega_1|$. Since $\lambda_1$ does not contain $\omega_1$ as a scattered subword, the number of occurrences $k_i$ of some variable $x_i$ in $\lambda_1$ must be less than in $\omega_1$, so the number of occurrences $k_j$ of some other variable $x_j$ must be greater  than the number of occurrences in $\omega_1$. Suppose $\mathcal{D}(\omega_2)$ is infinite. Then there exists a difference $\lambda_1'-\lambda_2' \in \mathcal{D}(\omega_2)$ with $|\lambda_1'| >|\lambda_1|$, and furthermore, neither $\lambda_1$ nor $\omega_1$ can occur as scattered subwords in $\lambda_1'$. Thus the number of occurrences $k_{i'}$ of another variable $x_{i'}$ must be less than in $\omega_1$, and so the number of occurrences $k_{j'}$ of another variable $x_{j'}$ must be greater than in $\omega_1$. This indicates that $\mathcal{D}(\omega_2)$ cannot be infinite: for some $l$, any homogeneous binomial difference $\gamma_1-\gamma_2 \in \mathcal{D}(\omega_2)$ such that  $|\gamma_1| > l$ must have a first word which contains as a scattered subword some word $\bar{\lambda}_1$ which previously occurred in a homogeneous binomial difference $\bar{\lambda}_1-\bar{\lambda}_2 \in \mathcal{D}(\omega_2)$ and is thus not minimal. The same argument, \emph{mutatis mutandis}, shows that $\mathcal{D}(\omega_1)$ is also finite. 

Consider, then, the set $\mathcal{D}(0)$. Let $\omega_1' - \omega_2' \in \mathcal{D}(0)$ be such that $|\omega_1'| \leq |\lambda_1|$ for any $\lambda_1-\lambda_2 \in \mathcal{D}(0)$. Note that $\omega_1' - \omega_2'$ must consist of words at least as long as $\omega_1$, and furthermore, in both $\omega_1'$ and $\omega_2'$ some variables $x_{k_1}$ and $x_{k_2}$ must occur less often than in $\omega_1$ and $\omega_2$ respectively. We may partition $\mathcal{D}(0)$ into sets $\mathcal{D}(\omega_1')$, $\mathcal{D}(\omega_2')$, and $\mathcal{D}(0')$ which form a partition of $\mathcal{D}(0)$. As above, these sets form a partition of $\mathcal{D}(0)$, and following the argument above, both $\mathcal{D}(\omega_1')$ and $\mathcal{D}(\omega_2')$ are finite. Consider then $\mathcal{D}(0')$, which must be infinite, and select a difference $\omega_1'' - \omega_2'' \in \mathcal{D}(0')$ such that $|\omega_1''|\leq |\lambda_1|$ for any $\lambda_1-\lambda_2 \in \mathcal{D}(0')$. Again $\omega_1'' - \omega_2''$ must consist of words at least as long as $\omega_1'$, and furthermore,  in both $\omega_1''$ and $\omega_2''$ some variables $x_{k'_1}$ and $x_{k'_2}$ must occur less often than in $\omega_1'$ and $\omega_2'$ respectively. Continuing this partitioning process \emph{ad infinitum} is impossible: for some $l$, any difference $\gamma_1 - \gamma_2$ such that $|\gamma_1| >l$ must contain the occurrence of some $\bar{\lambda}_i$, $i\in\{1,2\}$, which previously occurred in a homogeneous binomial difference in $\bar{\lambda}_1-\bar{\lambda}_2 \in \mathcal{D}(0)$ as a scattered subword. Thus $\mathcal{D}(0)$ cannot be infinite, and so $\mathcal{M}_A$ is finite and $I_A$ must be finitely generated.\qed
\end{pf}

We have the following corollaries. Corollary \ref{COR:Pairwise} gives a description of those homogeneous binomial differences in the commutator ideal. Note that necessity in Corollary \ref{COR:Pairwise} was proved in \citet{EH10}. 

\begin{cor}\label{COR:IAEqualsC} Let $A$ be a regular array with pairwise-coprime first column entries. Then $I_A = \mathcal{C}$, where $\mathcal{C}$ denotes the commutator ideal.
\end{cor}

\begin{pf} Since the entries in the first column of $A$ are pairwise-coprime, $\mathcal{M}_A$ is trivial. \qed 
\end{pf}

\begin{cor}\label{COR:Pairwise} Let $A$ be a regular array with pairwise-coprime first column entries, and suppose $\omega_1$, $\omega_2 \in \mathcal{B}$ with $|\omega_1| = |\omega_2| = l$. Then $\omega_1-\omega_2 \in I_A \iff \text{supp}(\omega_1)=\text{supp}(\omega_2)$ and $\#(x_i,\omega_1) = \#(x_i,\omega_2)$ for all $x_i \in \text{supp}(\omega_1) = \text{supp}(\omega_2)$. 
\end{cor}

\begin{pf} To prove sufficiency, let
\[
A=\left(\begin{array}{ccccc}
a_{1,0} & d_1a_{1,0} & \cdots & d_na_{1,0} & \cdots \\
a_{2,0} & d_1a_{2,0} & \cdots & d_na_{2,0} & \cdots \\
\vdots & \vdots & & \vdots & \\
a_{t,0} & d_1a_{t,0} & \cdots & d_na_{t,0} & \cdots 
\end{array}\right).
\]

Assume that  $\omega_1-\omega_2 \in I_A$. Then 
\[
\sigma(\omega_1) = \sigma(\omega_2) \Rightarrow \prod_{i = 0}^{l-1} a_{{u_i}i} = \prod_{i=0}^{l-1}a_{{v_i}i}.
\]

Expressing each weight as a multiple of a first-column entry and canceling the $d_i$'s common to each side of the equation gives
\begin{equation}\label{EQ:Pairwise}
\prod _{i = 0}^{l-1} a_{{u_i}0} = \prod _{i = 0}^{l-1} a_{{v_i}0}.
\end{equation}

Since the first column entries of $A$ are pairwise-coprime, equality can only hold in Equation (\ref{EQ:Pairwise}) when there exists a bijection between $\{a_{{u_i}0}\}_{i = 0}^{l-1}$ and $\{a_{{v_i}0}\}_{i = 0}^{l-1}$. Each $a_{{u_i}0}$ corresponds to an occurrence of the letter $x_{u_i}$ in $\omega_1$; thus, $\text{supp}(\omega_1) = \text{supp}(\omega_2)$, and because the weights are equal, $\#(x_i, \omega_1) = \#(x_i, \omega_2)$ for each $x_i \in \text{supp}(\omega_1)=\text{supp}(\omega_2)$. \qed
\end{pf}

In particular, when a regular array $A$ has pairwise-coprime first column entries, then the algebra $R/I_A$ on which it defines an order is in fact isomorphic to the commutative polynomial algebra $k[x_1, \cdots, x_t]$. The preceding construction can be viewed as an alternative way to specify an admissible length-dominant weight order on the monomials in $k[x_1, \dots, x_t]$. L. Robbiano has proven that any such order is a lexicographic product of weight orders \citep{LR86}.

More generally, in \citet{RG84}, R. Gilmer proved that monoid algebras of commutative monoids are precisely the homomorphic images of polynomial rings by ideals which are generated by pure binomial differences. Of course, the algebra $R/I_A$ is precisely such an algebra when $A$ is regular. Thus, we can view $R/I_A$ as a monoid algebra of a commutative monoid for which an admissible order on the basis of $R/I_A$ can be obtained. 

While Theorem \ref{THM:RegGen} shows how to decompose a homogeneous binomial difference over the set of commutators and $\mathcal{M}_A$ and Theorem \ref{THM:FinGen} demonstrates that $\mathcal{M}_A$ is finite, constructing $\mathcal{M}_A$ may present significant computational difficulties. Neglecting minimality, the process of directly identifying a homogeneous binomial difference of disjoint support consisting of words of length $l$ by calculating weights is easily seen to be equivalent to the Subset Product problem, which is known to be NP-complete \citep{MG79}. Furthermore, the above proofs do not give a bound on the lengths of the words that may occur in a minimal homogeneous binomial difference of disjoint support, suggesting that even if one is able to devise an algorithm to efficiently identify minimal homogeneous binomial differences of disjoint support, one could not terminate the algorithm and be satisfied that all the minimal homogeneous binomial differences of disjoint support had been enumerated.

\section{Weight Ideals Associated to Log-Linear Arrays}

Bijectively related to the family of linear arrays is the family of log-linear arrays. 

\begin{defn} An array $A=(a_{ij}) \in\mathcal{M}_{t \times \infty}\left(S_{>0}\right)$ is {\em log-linear\/} if the array $\log A = (\log{(a_{ij})})$ is linear. 
\end{defn}

\begin{exmp} The array 
\[
B=\left(\begin{array}{cccc}
e^2 & e^6 & e^{18} &\cdots  \\
e^3 & e^9 & e^{27} & \cdots \\
e^4 & e^{12} & e^{36} & \cdots 
\end{array}\right)
\]
is a log-linear array, because 
\[
\log B=\left(\begin{array}{cccc}
2 & 6 & 18 &\cdots  \\
3 & 9 & 27 & \cdots \\
4 & 12 & 36 & \cdots 
\end{array}\right)
\]
is a linear array (with slope $d=3$). 
\end{exmp}

Note that an array $A$ for which $\log A$ is regular but not linear need not be admissible. Because every log-linear array with slope $1$ is in fact a (constant) regular array, we will assume without comment in the remainder that any log-linear array considered has slope $d \neq1$. We point out also that the base of a log-linear array is immaterial in determining the order given by the array. To see this, let $A \in\mathcal{M}_{t \times \infty}\left(S_{>0}\right)$ be a given linear array with $i,j^{th}$ entry $d^ja_{i,0}$ and consider the arrays $B = (b_{ij})$ and $C = (c_{ij})$, where $b_{ij} = b^{d^ja_{i,0}}$ and $c_{ij} = c^{d^ja_{i,0}}$ for elements $b$, $c \in S_{>0}$ with $b \neq c$. Suppose that $\omega_1 = x_{u_0}x_{u_1}\dots x_{u_{l-1}}$ and $\omega_2 = x_{v_0}x_{v_1}\dots x_{v_{l-1}}$ are words in $X^*$ such that  $\omega_1 \succ_B \omega_2$. Then
\[
\sigma_B(\omega_1) > \sigma_B(\omega_2)
\]
which means that 
\[
\prod_{k=0}^{l-1}b_{u_k,k} > \prod_{i=0}^{l-1}b_{v_k,k}.
\]
Of course, this is equivalent to the inequality
\[
b^{\sum_{k=0}^{l-1}d^ka_{u_k}} > b^{\sum_{k=0}^{l-1}d^ka_{v_k}},
\]  
and replacing $b$ with $c$ does not change the direction of the inequality. Thus, we will typically assume without comment that the base of a log-linear array is $e$. As the preceding discussion indicates, the significant distinction between regular arrays and log-linear arrays is that when working with regular arrays, the weight associated to a word is calculated by multiplying the weights given to each variable in their respective positions, while when working with log-linear arrays, the weight associated to a word is calculated by adding the weights given to each variable in their respective positions. In particular, when working with a log-linear array $A$, the equation $\sum_{k=0}^{l-1}d^ka_{u_k} = \sum_{k=0}^{l-1}d^ka_{v_k}$ must be satisfied for a homogeneous binomial difference $\omega_1 - \omega_2 = x_{u_0}x_{u_1}\dots x_{u_{l-1}} - x_{v_0}x_{v_1}\dots x_{v_{l-1}}$ to be a member of $I_A$.

The structure of log-linear arrays is much less uniform than that of regular arrays. Theorem \ref{THM:Examples} is the main result of this section and is proved via the examples that follow. We will also demonstrate that log-linear arrays can be constructed to give orders on $R$ which are equivalent to the familiar left and right length-lexicographic orders.

\begin{thm}\label{THM:Examples} There exist log-linear arrays whose associated weight ideals are trivial, log-linear arrays whose associated weight ideals admit a finite generating set, and log-linear arrays whose associated weight ideals do not admit a finite generating set.
\end{thm}

It would be of interest to find necessary and sufficient conditions on a log-linear array $A$ such that $I_A$ is trivial, is nontrivial but admits a finite generating set, or is nontrivial and does not admit a finite generating set. 

The following two lemmas describe distinct arrays that define orders on $R$ which are equivalent to left and right length-lexicographic order respectively. The hypotheses on the array $A$ is sufficient to insure in each case that $I_A$ is trivial.  

\begin{lem}\label{LEM:LeftLenLex} Suppose the variables $x_1,\dots, x_t$ are ordered. Let $A \in \mathcal{M}_{t \times \infty}\left(S_{>0}\right)$ be log-linear with first column $A_{(0)} = [e^{a_{1,0}},\dots,e^{a_{t,0}}]^T$for which the values of the first-column entries of $A$ reflect the order given to $x_1,\dots,x_t$; that is, $x_{i_1} < x_{i_2}$ if and only if $a_{{i_1},0} < a_{{i_2},0}$ also. Let $\alpha$ and $\beta$ denote the minimum and maximum nonzero first column differences of $\log A$ respectively; that is, 
\[
\alpha = \min\{|a_{i,0} - a_{j,0}|: 1 \leq i, j\leq t, i \neq j\}, 
\]
and
\[
	\beta = \max\{|a_{i,0} - a_{j,0}|: 1 \leq i, j \leq t, i \neq j\}.
\]
 Let $d$ be the slope of $\log A$. If $d <1$ and $\alpha > d\beta/(1-d)$, then $I_A$ is trivial and the order given by $A$ is the left length-lexicographic order.
\end{lem} 

\begin{pf} Assume the hypotheses, and assume that $\omega_1 = x_{u_0}x_{u_1}\dots x_{u_{l-1}}$ and $\omega_2 = x_{v_0}x_{v_1}\dots x_{v_{l-1}}$ are two words of equal length in $R$ such that $\omega_1-\omega_2 \in I_A$. This implies that
\begin{equation}\label{EQ:First}
e^{a_{u_0,0}+da_{u_1,0}+\dots+d^{l-1}a_{u_{l-1},0}} = e^{a_{v_0,0}+da_{v_1,0}+\dots+d^{l-1}a_{v_{l-1},0}}, 
\end{equation}
and thus  
\[
a_{u_0,0}+da_{u_1,0}+\dots+d^{l-1}a_{u_{l-1},0} = a_{v_0,0}+da_{v_1,0}+\dots+d^{l-1}a_{v_{l-1},0}. 
\]
This in turn implies that
\begin{equation}\label{EQ:Second}
a_{u_0,0}-a_{v_0,0} = \sum_{i=1}^{l-1}d^i\left(a_{u_i,0}-a_{v_i,0}\right).
\end{equation}
Suppose now that the first letters of $\omega_1$ and $\omega_2$ differ. The largest in absolute value that the right-hand side of Equation \ref{EQ:Second} can be is when each difference $a_{u_i,0}-a_{v_i,0} = \beta$, so the right-hand side has an upper bound at $\beta(d-d^l)/(1-d)$. The smallest the left-hand side of Equation \ref{EQ:Second} can be in absolute value is when $a_{u_0,0} - a_{v_0,0}= \alpha$, and by hypothesis, $\alpha> d\beta/(1-d) > (d-d^l)\beta/(1-d)$ for all $l$. This contradicts the assumption that $\omega_1-\omega_2 \in I_A$, so in fact $I_A = \{0\}$. Furthermore, the difference $a_{u_0,0} - a_{v_0,0}$ is greater in absolute value than any possible subsequent sum and hence determines the order between $\omega_1$ and $\omega_2$.

Now, note that if the first $k$ letters of $\omega_1$ and $\omega_2$ are the same, then those first $k$ letters contribute the same expression to either side of Equation \ref{EQ:First}, and thus play no role in determining the order between $\omega_1$ and $\omega_2$. Thus, when the first $k$ letters are the same, we may determine the order between $\omega_1$ and $\omega_2$ by simply applying the above argument to the truncated words $x_{u_k}x_{u_{k+1}}\cdots x_{u_{l-1}}$ and $x_{v_k}x_{v_{k+1}}\cdots x_{v_{l-1}}$ to note that the order on $\omega_1$ and $\omega_2$ is determined solely by the order between $a_{u_k,0}$ and $a_{{v_k},0}$.  

Because the order on the first column entries of $\log A$ is equivalent to the order on the variables $x_1,\dots,x_t$ the order given by $A$ is thus the left length-lexicographic order. \qed  
\end{pf}

\begin{lem}\label{LEM:RightLenLex} Suppose the variables $x_1,\dots, x_t$ are ordered. Let $A \in \mathcal{M}_{t \times \infty}\left(S_{>0}\right)$ be log-linear with first column $A_{(0)} = [e^{a_{1,0}},\dots,e^{a_{t,0}}]^T$for which the values of the first column entries of $A$ reflect the order given to $x_1,\dots,x_t$; that is, $x_{i_1} < x_{i_2}$ if and only if $a_{{i_1},0} < a_{{i_2},0}$ also. Let $\alpha$ and $\beta$ denote the minimum and maximum nonzero first column differences of $\log A$ respectively; that is, 
\[
\alpha = \min\{|a_{i,0} - a_{j,0}|: 1 \leq i, j\leq t, i \neq j\}, 
\]
and
\[
	\beta = \max\{|a_{i,0} - a_{j,0}|: 1 \leq i, j \leq t, i \neq j\}.
\]
 If $d >1$ and $\alpha > \beta/(d-1)$, then $I_A$ is trivial and the order given by $A$ is the right length-lexicographic order.
\end{lem}

\begin{pf} The same argument as in Lemma \ref{LEM:LeftLenLex} applies, \textsl{mutatis mutandis}.
\end{pf}

Orders constructed via admissible arrays with trivial weight ideals are simply admissible length-dominant orders on $R$. A set of invariants that fully characterize the admissible orders that can be defined on a noncommutative $k$-algebra such as $R$ has not yet been described, though results in this direction have been obtained \citep{ES94, PF97}. It is possible that array-based admissible orders may be of use in defining such a set of invariants.

We turn our attention next to an example of a log-linear array $A$ for which $I_A$ is nontrivial and admits a finite generating set.

\begin{exmp}\label{EX:FinGen} Let $A$ be the log-linear array such that
\[
\log A=
\left(
\begin{array}{cccc}
2 & 4 & 8 & \dots \\
3 & 6 & 12 & \dots \\
4 & 8 & 16 & \dots \\
6 & 12 & 24 & \dots
\end{array}
\right).
\]

The weight ideal $I_A$ associated to $A$ is nontrivial and is finitely generated.
\end{exmp}

Clearly the weight ideal associated to $A$ is nontrivial; for example, $x_1x_2 - x_3x_1 \in I_A$. Interestingly, any homogeneous binomial difference of length $l>2$ in $I_A$ can be reduced in at most two steps to an algebraic consequence of homogeneous binomial differences consisting of words whose maximum length is $l-1$. It follows inductively that any homogeneous binomial difference of length $l$ can be reduced to an algebraic consequence of homogeneous binomial differences of length 2; that is, for this particular $A$, the homogeneous binomial differences of length 2 in fact generate $I_A$. The proof of this proposition is straightforward but relies on a lengthy case-by-case analysis, which is included as an appendix.

The next example demonstrates that there exists log-linear arrays with a nontrivial associated weight ideal which admits no finite generating set. 

\begin{exmp}\label{EX:InfGen} Let $A$ be the log-linear array given by
\[\log A =\left(\begin{array}{cccc}
2 & 4 & 8 & \dots \\
4 & 8 & 16 & \dots \\
7 & 14 & 28 & \dots 
\end{array}\right).
\]
The weight ideal $I_A$ associated to $A$ is nontrivial and does not admit a finite generating set. 
\end{exmp}

In the proof we will use the term \emph{factor} to indicate a subword in which the letters occur consecutively, in order to alleviate any potential confusion with scattered subwords, which some authors refer to simply as subwords.

\begin{pf} Consider the homogeneous binomial difference
\[
\omega_1-\omega_2 := x_2x_3^nx_2 - x_1(x_2x_3)^{(n-2)/2}x_1x_2x_3,
\]
where $n\geq4$ is an even integer and $l=n+2$ is the length of $\omega_1$ and $\omega_2$. We will show first that for any such $n$, the difference given above is a member of $I_A$, and then we will demonstrate that $\omega_1-\omega_2$ is not an algebraic consequence of the shorter length differences in $I_A$ and so must belong to any generating set for $I_A$. Since this holds for all $n \geq 4$, this will prove that $I_A$ does not admit a finite generating set.

To demonstrate that $x_2x_3^nx_2 - x_1(x_2x_3)^{(n-2)/2}x_1x_2x_3 \in I_A$ for any $n \geq 4$, let us calculate the difference of the weights associated to $\omega_1-\omega_2$ respectively. Fix an even integer $n\geq 4$. The difference in weights associated to any homogeneous binomial difference by $A$ can be expressed as a polynomial:
\[
	\Delta:= a_0 + da_1 +d^2a_2 + \dots +d^{l-1}a_{l-1},
\]
where $a_k$ is the difference of the first-column entries associated to the letter in position $k$ in $\omega_1$ and $\omega_2$ respectively. For the given difference, 
\[
\Delta_{\omega_1-\omega_2}=2+2\cdot3+2^2\cdot0+2^3\cdot3+2^4\cdot0+\dots+2^{l-4}\cdot0+2^{l-3}\cdot4+2^{l-2}\cdot3+2^{l-1}\cdot(-3).
\]
To show that $\Delta_{\omega_1-\omega_2}=0$ regardless of the value of $n$, it is easiest to work in binary. We have
\begin{multline*}
10+10\cdot11+1000\cdot11 +100000\cdot11+\cdots \\
+1\underbrace{0\dots0}_{l-5}\cdot11+1\underbrace{0\dots0}_{l-3}\cdot101+1\underbrace{0\dots0}_{l-2}\cdot11-1\underbrace{0\dots0}_{l-1}\cdot11.
\end{multline*}
Multiplying simplifies this to 
\[
10+110+11000+1100000+\dots+11\underbrace{0\dots0}_{l-3}+101\underbrace{0\dots0}_{l-3}+11\underbrace{0\dots0}_{l-2}-11\underbrace{0\dots0}_{l-1}.
\]
This expression is equal to 0:
\[
1\underbrace{0\dots0}_{l-3}+101\underbrace{0\dots0}_{l-3}+11\underbrace{0\dots0}_{l-2}-11\underbrace{0\dots0}_{l-1}
\]
and so 
\[
11\underbrace{0\dots0}_{l-2}+11\underbrace{0\dots0}_{l-2}-11\underbrace{0\dots0}_{l-1} =
\]
\[ 
11\underbrace{0\dots0}_{l-1}-11\underbrace{0\dots0}_{l-1} = 0.
\]

This demonstrates that $x_2x_3^nx_2 - x_1(x_2x_3)^{(n-2)/2}x_1x_2x_3 \in I_A$. To show that $x_2x_3^nx_2 - x_1(x_2x_3)^{(n-2)/2}x_1x_2x_3$ must be contained in any generating set for $I_A$, we will show that the word $x_2x_3^nx_2$ contains no factor that occurs as a word in a homogeneous binomial difference of shorter length in $I_A$. In particular, this implies that $x_2x_3^nx_2 - x_1(x_2x_3)^{(n-2)/2}x_1x_2x_3$ cannot be written as an algebraic consequence of strictly shorter length homogeneous differences in $I_A$. 

Consider the possible factors of the word $x_2x_3^nx_2$. For each $k$, $4\leq k\leq n$, we have factors $x_2x_3^k$, factors $x_3^kx_2$, and $x_3^k$. No factor that occurs in a homogeneous binomial difference in $I_A$ can begin with $x_3$, because of the parity of the weight that results, so we can rule out as possibilities any factors of the form $x_3^k$ and $x_3^kx_2$. The weight given the factor $x_2x_3^k$ will be greater than the weight assigned any other word of equal length except $x_3^k$, and it will not equal this weight. Thus, $x_2x_3^nx_2$ contains no factors that occur as a word in a homogeneous binomial difference in $I_A$. Because this holds for each $n\geq 4$, the difference $x_2x_3^nx_2 - x_1(x_2x_3)^{(n-2)/2}x_1x_2x_3$ must be included in any generating set for $I_A$, and thus any generating set for $I_A$ is infinite. \qed
\end{pf}

\appendix

\section{Proof that the Array in Example \ref{EX:FinGen} is Finitely Generated}
\begin{pf} Let us first list the the homogeneous binomial differences of length two that occur in $I_A$. They are $x_1x_2-x_3x_1$, $x_1x_3-x_3x_2$, $x_1x_3-x_4x_1$, $x_3x_2-x_4x_1$, $x_3x_3-x_4x_2$, and $x_1x_4-x_4x_3$. Given a homogeneous binomial difference consisting of words of length $l$ in $I_A$, we will show that it either decomposes over the homogeneous binomial differences in $I_A$ of length two or fails to belong to $I_A$. Proof of the latter claim requires us to consider the existence of solutions to the polynomial
\[
	a_0+da_1+\cdots+d^{l-1}a_{l-1} = 0  
\]
which corresponds to a given homogenous binomial difference in $I_A$. Any solution to this equation is an element of the set $\coprod_{i=0}^{l-1}A_{d} = \{(a_0,\dots,a_{l-1})| a_i = a_{j,0} - a_{k,0}\text{,}1\leq j,k\leq 4 \}$. 

To simplify exposition, we use the notation of rewriting relations on words in the free monoid $X^*:=\langle x_1,\dots,x_t\rangle$. We define the following rewriting relation: for $\omega_1$, $\omega_2 \in X^*$, we write $\omega_1\overset{*}\longleftrightarrow \omega_2$ to denote that $\omega_1-\omega_2 \in I_A$. A particular chain of rewritings $\omega_1\overset{*}\longleftrightarrow\lambda_1\overset{*}\longleftrightarrow \cdots\overset{*}\longleftrightarrow\lambda_n\overset{*}\longleftrightarrow\omega_2$ corresponds to a unique decomposition of $\omega_1-\omega_2$ in $R$ \citep{MR97}, though this decomposition need not be over homogeneous binomial differences consisting of words of lesser length. However, a chain of rewritings $\omega_1\overset{*}\longleftrightarrow\lambda$ where at least the first letter of $\lambda$ is the same as the first letter of $\omega_2$ does correspond to a unique decomposition of $\omega_1-\omega_2$ over the set of homogeneous binomial differences in $I_A$ whose words are of lesser length. Rather than calculate the decomposition precisely, we will rewriting to indicate that a decomposition is possible. 

Organizing by weight, the length two homogeneous binomial differences in $I_A$ give rise to the following rewriting relations:
\[
x_1x_2 \overset{*}\longleftrightarrow x_3x_1\text{, \hspace{0.25in}}x_1x_3\overset{*}\longleftrightarrow x_3x_2\overset{*}\longleftrightarrow x_4x_1,
\]

\[
x_3x_3\overset{*}\longleftrightarrow x_4x_2\text{,\hspace{0.25in}} x_1x_4 \overset{*}\longleftrightarrow x_4x_3.
\]

Suppose that $\omega_1 - \omega_2 \in I_A$, with $|\omega_1| = |\omega_2| = l$. Let $\omega_1 = x_{u_0}\dots x_{u_{l-1}}$ and $\omega_2 = x_{v_0}\dots x_{v_{l-1}}$. By parity, if either $\omega_1$ or $\omega_2$ start with $x_2$, then so must the other, in which case $\omega_1 - \omega_2$ is immediately reducible to a homogeneous binomial difference of length $l-1$. Thus we may assume without loss of generality that neither $\omega_1$ nor $\omega_2$ start with $x_2$.

Assume next that $\omega_1$ begins with $x_1$. We need to consider the cases when $\omega_2$ starts with $x_3$ or $x_4$, as well as the case when $\omega_1$ begins with $x_3$ and $\omega_2$ begins with $x_4$. All other cases will then be captured by symmetry.

If $\omega_1 - \omega_2$ is not immediately reducible, then $\omega_2$ must begin with either $x_3$ or $x_4$.

\noindent\textbf{Case 1} : $\omega_2$ begins with $x_3$.

\textbf{Subcase 1.1} : $\omega_2$ begins with $x_3x_1$.

$x_3x_1 \overset{*}\longleftrightarrow x_1x_2$, so $\omega _1-\omega_2$ is reducible after this single rewrite.

\textbf{Subcase 1.2} : $\omega_2$ begins with $x_3x_2$.

$x_3x_2 \overset{*}\longleftrightarrow x_1x_3$, so $\omega_1-\omega_2$ is reducible after this single rewrite. 

\textbf{Subcase 1.3} : $\omega_2$ begins with $x_3x_3$.

In this case, reduction with a single rewrite of $\omega_2$ is not always possible. If the second letter of $\omega_1$ is $x_2$ or $x_3$, then we can rewrite $\omega_1$ to reduce $\omega_1 - \omega_2$. Suppose then that the second letter of $\omega_1$ is $x_4$. Rewrite as follows:
\[
x_3x_3 \overset{*}\longleftrightarrow x_4x_2\text{ and }x_1x_4 \overset{*}\longleftrightarrow x_4x_3.
\]
Now first letters agree and $\omega_1-\omega_2$ is reducible. Finally, suppose that the second letter of $\omega_1$ is $x_1$. Then $\omega_1-\omega_2$ is of the form
\[
x_1x_1x_{u_2}\dots x_{u_{l-1}} - x_3x_3x_{v_2}\dots x_{v_{l-1}}.
\]

Since $\omega_1-\omega_2 \in I_A$, this gives rise to the equation $a_0+2a_1+\cdots+2^{l-1}a_{l-1}=0$, where $a_0$ is a difference of first-column entries of $I_A$ determined by the letters of $\omega_1$ and $\omega_2$. In particular, $a_0  = -2$ and $a_1 = -2$, so
\[
-2-4+2^2a_2+\cdots+2^{l-1}a_{l-1}=0,  
\]
or equivalently,
\begin{equation}\label{EQ:Impos}
2a_2+\cdots+2^{l-2}a_{l-1}=3,  
\end{equation}
but (\ref{EQ:Impos}) does not have a solution over $\coprod_{i=0}^{l-1}A_{d}$.

\textbf{Subcase 1.4} : $\omega_2$ begins with $x_3x_4$.

Again, if the second letter of $\omega_1$ is $x_2$ or $x_3$, then $\omega_1-\omega_2$ is immediately reducible. Suppose the second letter of $\omega_1$ is $x_1$. Then $\omega_1-\omega_2$ is of the form
\[
x_1x_1x_{u_2}\dots x_{u_{l-1}} - x_3x_4x_{v_2}\dots x_{v_{l-1}}.
\]

Since $\omega_1-\omega_2 \in I_A$, this gives rise to the equation $a_0+2a_1+\cdots+2^{l-1}a_{l-1}=0$, where each $a_i$ is a difference of first column entries of $A$ determined by the letters of $\omega_1$ and $\omega_2$. In particular, $a_0  = -2$ and $a_1 = -4$, so
\[
-2-8+2^2a_2+\cdots+2^{l-1}a_{l-1}=0,  
\]
or equivalently,
\begin{equation}\label{EQ:Impos2}
2a_2+\cdots+2^{l-2}a_{l-1}=5,  
\end{equation}
but (\ref{EQ:Impos2}) does not have a solution over $\coprod_{i=0}^{l-1}A_{d}$.

Finally, if the second letter of $\omega_1$ is $x_4$, then $\omega_1-\omega_2$ is of the form
\[
x_1x_4x_{u_2}\dots x_{u_{l-1}} - x_3x_4x_{v_2}\dots x_{v_{l-1}}.
\]
Since $\omega_1-\omega_2 \in I_A$, this gives rise to the equation $a_0+2a_1+\cdots+2^{l-1}a_{l-1}=0$, where each $a_i$ is a difference of first column entries of $A$ determined by the letters of $\omega_1$ and $\omega_2$. In particular, $a_0  = -2$ and $a_1 = 0$, so
\[
-2-0+2^2a_2+\cdots+2^{l-1}a_{l-1}=0,  
\]
or equivalently,
\begin{equation}\label{EQ:Impos3}
2a_2+\cdots+2^{l-2}a_{l-1}=1,  
\end{equation}
but (\ref{EQ:Impos3}) does not have a solution over $\coprod_{i=0}^{l-1}A_{d}$.

\noindent\textbf{Case 2}: $\omega_1$ begins with $x_1$ and $\omega_2$ begins with $x_4$.

\textbf{Subcase 2.1}: $\omega_1$ begins with $x_1$ and $\omega_2$ begins with $x_4x_1$.

$x_4x_1 \overset{*}\longleftrightarrow x_1x_3$, so $\omega_1-\omega_2$ is reducible after this single rewrite. 

\textbf{Subcase 2.2}: $\omega_1$ begins with $x_1$ and $\omega_2$ begins with $x_4x_2$.

If $\omega_1$ begins with $x_1x_3$, then we can rewrite $x_1x_3 \overset{*}\longleftrightarrow x_4x_1$ and immediately reduce $\omega_1-\omega_2$. Similarly, if $\omega_1$ begins with $x_1x_4$, we can immediately rewrite $x_1x_4 \overset{*}\longleftrightarrow x_4x_3$ to reduce $\omega_1-\omega_2$. If $\omega_1$ begins with $x_1x_2$, rewrite $x_1x_2 \overset{*}\longleftrightarrow x_3x_1$ and $x_4x_2 \overset{*}\longleftrightarrow x_3x_3$ to reduce $\omega_1-\omega_2$. The only remaining possibility is that $\omega_1$ begins with $x_1x_1$ and $\omega_2$ begins with $x_4x_2$, but the corresponding polynomial shows that no such difference that begins with these letters can belong to $I_A$:
\[
-4+2\cdot(-1)+2^2a_2+\cdots+2^{l-1}a_{l-1} = 0
\]
reduces to 
\[
-3+2a_2+\cdots+2^{l-2}a_{l-1} = 0,
\]
and this equation has no solution over $\coprod_{i=0}^{l-1}A_{d}$. 

\textbf{Subcase 2.3}: $\omega_1$ begins with $x_1$ and $\omega_2$ begins with $x_4x_3$.

$x_4x_3 \overset{*}\longleftrightarrow x_1x_4$, so $\omega_1-\omega_2$ is reducible after this single rewrite.

\textbf{Subcase 2.4}: $\omega_1$ begins with $x_1$ and $\omega_2$ begins with $x_4x_4$.

\hspace{0.25in}\textbf{Subsubcase 2.4.1}: $\omega_1$ begins with $x_1x_1$.

Consider the polynomial equation $a_0+2a_1+\cdots +2^{l-1}a_{l-1}=0$ corresponding to this difference. The choice of first letters for $\omega_1$ and $\omega_2$ determine $a_0 = -4$ and $a_1 = -4$. Factoring out 4 from the resulting equation gives
\[
-3+a_2+2a_3+\cdots+2^{l-3}a_{l-1} = 0.
\]
Each term after the second in the expression on the left--hand side of the above equation is congruent to 0 mod 2, thus this equation has a solution only if $a_2 = \pm 3$ or $\pm 1$ (each $a_i \in\{\pm1,\pm2,\pm3\pm4\}$. There are thus four possible subcases to consider. If $a_2 = -3$, then the first three letters of $\omega_1$ are $x_1x_1x_2$, and $x_2x_1x_2 \overset{*}\longleftrightarrow x_1x_3x_1 \overset{*}\longleftrightarrow x_4x_1x_1$, so $\omega_1-\omega_2$ is reducible. If $a_2 = 3$, the first three letters of $\omega_1$ are $x_1x_1x_4$, and $x_1x_1x_4 \overset{*}\longleftrightarrow x_1x_4x_3 \overset{*}\longleftrightarrow x_4x_4x_3$, so $\omega_1-\omega_2$ is reducible. If $a_2 = -1$ then either $\omega_1$ begins with $x_1x_1x_1$ or $x_1x_1x_2$. In the first case, $\omega_2$ must therefore begin with $x_4x_4x_2$, and $x_4x_4x_2 \overset{*}\longleftrightarrow x_4x_3x_3 \overset{*}\longleftrightarrow x_1x_4x_3$, so the difference is reducible, and in the second case, $x_1x_1x_2 \overset{*}\longleftrightarrow x_1x_3x_1 \overset{*}\longleftrightarrow x_4x_1x_1$, so the difference is reducible. If $a_2 = 1$, then either $\omega_1$ begins with $x_1x_1x_2$ or $\omega_1$ begins with $x_1x_1x_3$. The first case has already been addressed, and the second case gives rise to a reducible instance of $\omega_1-\omega_2$, for $x_1x_1x_3 \overset{*}\longleftrightarrow x_1x_3x_2 \overset{*}\longleftrightarrow x_4x_1x_2$.

\hspace{0.25in}\textbf{Subsubcase 2.4.2}: $\omega_1$ begins with $x_1x_2$.

The polynomial equation corresponding to this difference is
\[
	-4+2(-3)+\cdots+2^{l-1}a_{l-1} = 0,
\]
and this equation has no solution over $\coprod_{i=0}^{l-1}A_{d}$. 

\hspace{0.25in}\textbf{Subsubcase 2.4.3}: $\omega_1$ begins with $x_1x_3$. 

$x_1x_3\overset{*}\longleftrightarrow x_4x_3$, so $\omega_1-\omega_2$ is reducible after this single rewrite.

\noindent\textbf{Case 3}: $\omega_1$ begins with $x_3$ and $\omega_2$ begins with $x_4$.

\textbf{Subcase 3.1}: $\omega_1$ begins with $x_3$ and $\omega_2$ begins with $x_4x_1$.

$x_4x_1 \overset{*}\longleftrightarrow x_3x_2$, so $\omega_1-\omega_2$ is reducible after this single rewrite.

\textbf{Subcase 3.2}: $\omega_1$ begins with $x_3$ and $\omega_2$ begins with $x_4x_2$.

$x_4x_2 \overset{*}\longleftrightarrow x_3x_3$, so $\omega_1-\omega_2$ is reducible after this single rewrite.

\textbf{Subcase 3.3}: $\omega_1$ begins with $x_3$ and $\omega_2$ begins with $x_4x_3$.

There are a number of cases to consider. If $\omega_1$ begins with $x_3x_1$, since $x_3x_1 \overset{*}\longleftrightarrow x_1x_2$ and $x_4x_3 \overset{*}\longleftrightarrow x_1x_4$, the difference is reducible. If $\omega_1$ begins with $x_3x_2$, we can rewrite $x_3x_2\overset{*}\longleftrightarrow x_1x_3$ and $x_4x_3 \overset{*}\longleftrightarrow x_1x_4$ to reduce $\omega_1-\omega_2$. If $\omega_1$ begins with $x_3x_3$, we may rewrite $x_3x_3 \overset{*}\longleftrightarrow x_4x_2$ to immediately reduce $\omega_1-\omega_2$. It remains to consider the case when $\omega_1$ begins with $x_3x_4$ and $\omega_2$ begins with $x_4x_3$. The corresponding polynomial shows that a homogeneous binomial difference that starts with these letters cannot occur in $I_A$:
\[
-2+2\cdot(2)+2^2a_2+\cdots+2^{l-1}a_{l-1} = 0
\] 
implies 
\[
1+2a_2+\cdots+2^{l-2}a_{l-1} = 0,
\]
and this equation does not have a solution over $\coprod_{i=0}^{l-1}A_{d}$.

\textbf{Subcase 3.4}: $\omega_1$ begins with $x_3$ and $\omega_2$ begins with $x_4x_4$.

The corresponding polynomial equation for a difference with these starting letters is 
\[
-2+2a_1+\cdots+2^{l-1}a_{l-1} = 0,
\]
and we can factor out the common 2 to obtain
\[
-1+a_1+\cdots+2^{l-1}a_{l-1}=0.
\]
In order for this equation to have a solution, $a_1 \in \{\pm1,\pm3\}$, but because the second letter of $\omega_2$ is $x_4$, $a_1 \neq\pm1$ and $a_1 \neq 3$, so $a_1 = -3$. Thus $\omega_1$ begins with $x_3x_2$, and because $x_3x_2 \overset{*}\longleftrightarrow x_4x_1$, $\omega_1-\omega_2$ is reducible. \qed 
\end{pf}


\end{document}